\documentclass[12pt]{amsart}
\usepackage{amssymb,amsmath}
\headsep=1cm \numberwithin{equation}{section}
\theoremstyle{definition}

\theoremstyle{remark}

\newcommand{\Ext}{\operatorname{Ext}}
\newcommand{\ann}{\operatorname{Ann}}
\newcommand{\Spec}{\operatorname{Spec}}

\newcommand{\Ass}{\operatorname{Ass}}

\newcommand{\Assh}{\operatorname{Assh}}
\newcommand{\Att}{\operatorname{Att}}

\newcommand{\Supp}{\operatorname{Supp}}

\newcommand{\fm}{\mathfrak{m}}
\newcommand{\fp}{\mathfrak{p}}
\newcommand{\fq}{\mathfrak{q}}
\newcommand{\fb}{\mathfrak{b}}
\newcommand{\fa}{\mathfrak{a}}

\newcommand{\LH}{\operatorname{H}}
\newcommand{\cd}{\operatorname{cd}}

\begin{document}

\author[M. T. Dibaei and S. Yassemi]{Mohammad T. Dibaei and
 Siamak Yassemi}
\title[Attached primes of the top local cohomology modules $\cdots$ (II)]
 { Attached primes of the top local cohomology modules with respect to an ideal (II)}


\begin{abstract} For a finitely generated module $M$, over a commutative
Noetherian local ring $(R, \fm)$, it is shown that there exist
only a finite number of non--isomorphic top local cohomology
modules $\LH _{\fa}^{\dim (M)}(M)$, for all ideals $\fa$ of $R$.
We present a reduced secondary representation for the top local
cohomology modules with respect to an ideal. It is also shown that
for a given integer $r\geq 0$, if $\LH _{\fa}^{r}(R/\fp)$ is zero
for all $\fp$ in $\Supp (M)$, then $\LH _{\fa}^{i}(M)= 0$ for all
$i\geq r$.

\end{abstract}

\maketitle

\section*{Introduction}
Throughout $(R, \fm)$ denote a commutative Noetherian local ring
and $\fa$, $\fb$ are proper ideals of $R$. Assume that $M$ is a
finite (i.e. finitely generated) $R$--module of dimension $n$.
Recall that for an $R$--module $T$, a prime ideal $\fp$ of $R$ is
 said to be an attached prime of $T$ if $\fp=\ann (T/U)$ for
 some submodule $U$ of $T$. We denote the set of attached primes
 of $T$ by $\Att (T)$. This definition agrees with the usual
 definition of attached primes if $T$ has a secondary
 representation [M, Theorem 2.5].

 In this paper, we are interested in the structure of
${\LH}_{\fa}^{n}(M)$, which is the $n$--th right derived functor
of $\Gamma_{\fa}(M)= \underset{i>0}\bigcup 0\underset{M}: \fa
^{i}$ on $M$.  In [MS] Macdonald and Sharp studied the top local
cohomology module with respect to the  maximal ideal and showed
that
  ${\LH}_{\fm}^{n}(M)$ is an Artinian module and ${\Att}({\LH}_{\fm}^{n}(M))=
  \{\fp \in {\Ass}(M): {\dim}(R/\fp) = n\}$
    (the right hand set is denoted by ${\Assh}(M)$). In [S] Sharp proved that
    ${\LH}_{\fm}^{n}(M)$ is isomorphic to the last term of the Cousin Complex of
  $M$, $C_{R}(M)$, with respect to the height filtration of $M$.  \\
  \indent It is well--known that for any integer $i$ there is an exact sequence
  $$\LH_{\fm}^{i}(M)
  \longrightarrow \LH_{\fa}^{i}(M)\longrightarrow
  \underset{\underset{t}\longrightarrow}\lim \Ext _{R}^{i}(\fm ^{t}/\fa ^{t}, M)$$
  relating to the local cohomology modules of $M$ with respect to $\fa$ and $\fm$, cf. [G].
  By Hartshorne's  result [H, Page 417], we know that the right hand side module in the above exact sequence is zero.
  That means there exists an epimorphism $\LH_{\fm}^{n}(M) \longrightarrow
  \LH_{\fa}^{n}(M)$ ( in [DS, Theorem 3.2], Divaavi-Aazar and Schenzel
  calculated the kernel of this epimorphism). In [DY], we
  investigated the attached primes of ${\LH}_{\fa}^{n}(M)$, and
  showed  that
  $${\Att}({\LH}_{\fa}^{n}(M))=\{ \fp\in {\Supp}(M)| \LH _{\fa}^{n}(R/\fp) \neq 0 \}.$$

This paper is divided into two sections. In section 1, we continue
our investigation of the top local cohomology module
$\LH_{\fa}^{n}(M) $. We first show that if $R$ is complete then
${\LH}_{\fa}^{n}(M)$ is isomorphic to ${\LH}_{\fm }^{n}(K)$ for
some homomorphic image $K$ of $M$ (see Theorem 1.3). Then we give
a generalization of Hartshorne's result by showing that if
$\fa\subseteq \fb$, then there exists an epimorphism
$\LH_{\fb}^{n}(M) \longrightarrow \LH_{\fa}^{n}(M)$ (see Theorem
1.5).

 It is also shown that,
 if ${\Att}({\LH}_{\fa}^{n}(M))= {\Att}({\LH}_{\fb}^{n}(M))$ then
${\LH}_{\fa}^{n}(M)\cong {\LH}_{\fb}^{n}(M)$ (see Proposition
1.7). This result implies that the number of non--isomorphic top
local cohomology modules of $M$ is less than or equal to
$2^{|\Assh(M)|}$. We close this section by presenting a reduced
secondary representation of $\LH^n_\fa(M)$ (see Theorem 1.10).

It is not always the case  that ${\LH}_{\fa}^{n}(M)\neq 0$, so the
supremum of all integers $i$ for which ${\LH}_{\fa}^{i}(M)\neq 0$,
denoted by $\cd (\fa, M)$, is called the cohomological dimension
of $M$ with respect to $\fa$. It seems that it is not so easy to
know the structure of ${\LH}_{\fa}^{\cd (\fa, M)}(M)$. In [H,
Proposition 2.3], Hartshorne showed that for an integer $r$, if
$\LH _{\fa}^{r}(N)= 0$ for all $R$--modules $N$, then the
corresponding condition holds for all $j\geq r$. In section 2, for
a non--negative integer $r$, we generalize Hartshorne's result by
showing that $\cd (\fa, M)< r$ if and only if
${\LH}_{\fa}^{r}(R/\fp) = 0$ for all $\fp\in\Supp (M)$ (see
Theorem 2.1). Then we specify a subset of $\Att ({\LH}_{\fa}^{\cd
(\fa, M)}(M))$ by showing that (see Theorem 2.4) $$\{ \fp\in
\Ass(M) : \cd (\fa, R/\fp)= \dim (R/\fp)= \cd (\fa,
 M)\}\subseteq \Att (\LH_{\fa}^{\cd (\fa,
 M)}(M)).$$

\section{Top local cohomology modules}

In this section we investigate  the top local cohomology module
$\LH_{\fa}^{n}(M) $. For a module T over a commutative ring $S$
and $X\subseteq \Ass(T)$, there exists a submodule $T'$ of $T$
such that $\Ass(T')=\Ass(T)\setminus X$ and $\Ass (T/T')= X$ (see
[B, p. 263, Proposition 4]). We use this fact frequently
throughout the paper without further comment. Recall the following
result from [DY, Theorem A].
\subsection{Theorem }
{\it We have}\\
$$\Att (\LH_{\fa}^{n}(M)) = \{\fp\in \Supp (M) : \cd (\fa, R/\fp)= n \}.$$

We first prove the following result which is essential for our
 approach.
 \subsection{Lemma }
 {\it Assume that $\LH_{\fa}^{n}(M)\neq 0$. Then there exists a homomorphic image $K$ of
 $M$, such that \\
 (i) $\dim (K) = n$,\\
 (ii) $K$ has no non--zero submodule of dimension less than $n$,\\
 (iii) ${\Ass}(K) = \{\fp \in \Ass (M) : \cd (\fa, R/\fp) = n\}$,\\
 (iv) ${\LH}_{\fa}^{n}(K) \cong {\LH}_{\fa}^{n}(M)$,\\
 (v) ${\Ass}(K)={\Att}({\LH}_{\fa}^{n}(K))$} .\\

 {\it Proof}. By Theorem 1.1, ${\Att}({\LH}
 _{\fa}^{n}(M))=\{ \fp\in {\Supp}(M): {\cd}(\fa, R/\fp)=n\}$. Note
 that if $\fp\in {\Att}({\LH}
 _{\fa}^{n}(M))$, then we have
 $$n={\cd}(\fa, R/\fp)\leq {\dim}(R/\fp)\leq {\dim}(M)=n,$$
  and so ${\Att}({\LH} _{\fa}^{n}(M))\subseteq {\Ass}(M)$.
  Therefore there is a submodule $N$ of $M$ such that  ${\Ass}(N)={\Ass}(M)\setminus {\Att}
  ({\LH} _{\fa}^{n}(M))$ and ${\Ass}(M/N)=
  {\Att}({\LH} _{\fa}^{n}(M))$ .  Now, consider the exact sequence
  ${\LH}_{\fa}^{n}(N)\longrightarrow {\LH}_{\fa}^{n}(M)\longrightarrow
  {\LH}_{\fa}^{n}(M/N)\longrightarrow 0$. We claim that
   ${\LH}_{\fa}^{n}(N)= 0$. In the other case,  by Theorem 1.1, there
   exists
  $\fp\in {\Att}({\LH} _{\fa}^{n}(N))$ (and hence $\fp\in\Ass (N)$) such that ${\cd}(\fa, R/\fp)=n$.
  Therefore $\fp\in \Ass (M)$ and hence $\fp\in \Att(\LH _{\fa}^{n}(M))$, that is a
contradiction.
  Thus ${\LH}_{\fa}^{n}(M)\cong
  {\LH}_{\fa}^{n}(M/N)$. Now the assertions holds for $K:= M/N$. \hfill $\Box$\\

  Now we are ready to present one of our main results. In the
  following $\widehat{L}$ indicates the completion of an $R$-module
  $L$ with respect to the $\fm$-adic topology.\\

\subsection{Theorem }
{\it Assume that ${\LH} _{\fa}^{n}(M)\neq 0$. Then there exists a
homomorphic image $K$ of $\widehat{M}$ with ${\Ass}_{ \widehat{R
}}(K)={\Att}_{ \widehat{R}}({\LH} _{\fa}^{n}(M))$ and ${\LH}
_{\fa}^{n}(M)\cong {\LH} _{\fm \widehat{R}}^{n}(K)$ as $ \widehat{R}$--modules.}\\

{\it Proof.}  As
  ${\LH}_{\fa}^{n}(M)$ is an Artinian $R$--module, it has the
  structure of $\widehat{R}$--module and that
  ${\LH}_{\fa}^{n}(M)\cong {\LH}_{\fa
  \widehat{R}}^{n}(\widehat{M})$. By applying Lemma 1.2, there
  exists a finite $ \widehat{R}$--module $K$ which is a homomorphic image of
   $ \widehat{M}$ such that ${\Ass}
  _{\widehat{R}}(K)={\Att} _{\widehat{R}}({\LH}_{{\fa} \widehat{R}}^{n}(\widehat{M}))$ and
${\LH}_{\fa \widehat{R}}^{n}(K)\cong {\LH}_{\fa
\widehat{R}}^{n}(\widehat{M})$.\\
\indent If $\fq \in {\Ass} _{\widehat{R}}(K)$, then ${\cd}(\fa
\widehat{R}, \widehat{R}/\fq)=n={\dim} (\widehat{R}/\fq)$ and, by
Lichtenbaum--Hartshorne Vanishing Theorem of local cohomology,
$\fa \widehat{R}+\fq$ is an $\fm \widehat{R}$--primary ideal.
Therefore, ${\ann}_{ \widehat{R}}(K)+ \fa \widehat{R}$ is an $\fm
\widehat{R}$--primary ideal too. Hence we have
$${\LH} _{\fa \widehat{R}}^{n}(K)\cong {\LH} _{\fa
\widehat{R}+
{\ann}_{\widehat{R}}(K)/{\ann}_{\widehat{R}}(K)}^{n}(K)\cong {\LH}
_{\fm \widehat{R}/{\ann}_{\widehat{R}}(K)}^{n}(K).$$  Thus , by
using the Independence Theorem of local cohomology (cf. [BS,
4.2.1]), ${\LH}
_{\fa}^{n}(M)\cong {\LH} _{\fm \widehat{R}}^{n}(K)$. \hfill $\Box$\\

It is well known that for a ring $R$ with positive dimension if
${\LH}^{\dim (R)}_\fa(R)$ is non-zero then it is not finite, cf.
[BS, Exercise 8.2.6 (ii)]. As a simple corollary of Theorem 1.3 we
give the next result.

\subsection{Corollary }
{\it If $n>0$ and ${\LH} _{\fa}^{n}(M)\neq 0$, then ${\LH}
_{\fa}^{n}(M)$ is not finitely generated $R$--module.}\\

{\it Proof.} There exists a finite $\widehat{R}$-module $K$ such
that $\LH^n_\fa(M)\cong\LH^n_{\fm\widehat{R}}(K)$. Now the
assertion holds by the fact that $\LH^n_{\fm\widehat{R}}(K)$ is not finitely generated. \hfill $\Box$\\

As mentioned in the introduction, by using the results of
Grothendieck [G] and Hartshorne [H, page 417], there is an
epimorphism $\LH _{\fm}^{n}(M) \longrightarrow \LH _{\fa}^{n}(M)$.
In the following result we generalize it for two ideals
$\fa\subseteq \fb$.
\subsection{Theorem}
{\it If $\fa\subseteq \fb$, then there is an epimorphism $\LH
_{\fb}^{n}(M) \longrightarrow \LH _{\fa}^{n}(M)$}.\\

{\it Proof}. We may assume that $\fa\neq \fb$ and choose $x\in
\fb\setminus \fa$. By [BS, Proposition 8.1.2], there is an exact
sequence
$$\LH _{\fa +xR}^{n}(M)\longrightarrow \LH
_{\fa}^{n}(M)\longrightarrow \LH _{\fa}^{n}(M_{x}),$$ where
$M_{x}$ is the localization of $M$ at $\{x^{i}: i\geq 0\}$. Note
that $\LH _{\fa}^{n}(M_{x})\cong \LH _{\fa R_{x}}^{n}(M_{x})$ and
$\dim _{R_{x}}(M_{x})<n$, so that $\LH _{\fa R_{x}}^{n}(M_{x})=
0$. Thus there is an epimorphism $\LH _{\fa
+xR}^{n}(M)\longrightarrow \LH _{\fa}^{n}(M)$. Now the assertion
follows by assuming ${\fb} = {\fa} +(x_{1},\cdots,
x_{r})$ and applying the argument for finite steps. \hfill $\Box$\\

For ideals $\fa\subseteq\fb$ it follows from Theorem 1.5 that
${\Att}({\LH}_{\fa}^{n}(M)) \subseteq {\Att}({\LH}_{\fb}^{n}(M))$.
In the following result we give another proof (in some sense) of
this fact.

\subsection{Proposition}
{\it Let $\fa$ and $ \fb$ be ideals of $R$ such that $\fa
\subseteq \fb$. Then ${\Att}({\LH}_{\fa}^{n}(M)) \subseteq
{\Att}({\LH}_{\fb}^{n}(M))$. In particular, if ${\LH}_
{\fa}^{n}(M) \neq 0$, then ${\LH}_{\fb}^{n}(M) \neq 0$.}\\

{\it Proof}. We may assume that ${\LH}_ {\fa}^{n}(M) \neq 0$ and
$\fp\in {\Att}({\LH}_{\fa}^{n}(M))$. As mentioned in the
introduction, ${\cd}(\fa, R/\fp)=n$. Set $S=R/\ann M$. By [DY,
Theorem A and Corollary 4], there is $\fq \in {\Assh}
(\widehat{S})$ such that $\fp=\fq\cap R$ and ${\dim}(\widehat{S}/
\fa \widehat{S} +\fq) =0$. This implies that
${\dim}(\widehat{S}/\fb \widehat{S} +\fq) =0$. Again, by [DY,
Corollary 4 and Theorem A], $\fp\in
{\Att}({\LH}_{\fb}^{n}(M))$.\hfill $\Box$\\

In the following result it is shown that if $R$ is complete, then
the existence of an epimorphism $\LH^n_\fb(M)\to\LH^n_\fa(M)$
still
valid under weaker condition. More precisely:\\

\subsection{Theorem}
{\it Assume that $R$ is complete and that
$\Att(\LH^n_\fa(M))\subseteq\Att(\LH^n_\fb(M))$. Then there exists
an
epimorphism $\LH^n_\fb(M)\to\LH^n_\fa(M)$.}\\

{\it Proof}. Choose a submodule $N$ of $M$ with $\Ass (N)=\Ass
(M)\setminus\Att \LH^n_\fb(M)$ and $\Ass
(M/N)=\Att(\LH^n_\fb(M))$. As in the proof of Lemma 1.2,
$\LH^n_\fb(N)=0$. Hence we have an isomorphism
$\LH^n_\fb(M)\cong\LH^n_\fb(M/N)$. Note that any element
$\fp\in\Ass (M/N)$ satisfies $\cd(\fb, R/\fp)=n$, so that, by
Lichtenbaum--Hartshorne Vanishing Theorem, $\fb+\fp$ is
$\fm$-primary. This shows that $\fb+\ann(M/N)$ is an $\fm$-primary
ideal and thus
$$\LH^n_\fb(M/N)\cong\LH^n_{\fb+\ann(M/N)}(M/N)\cong\LH^n_\fm(M/N).$$
On the other hand, in the exact sequence
$\LH^n_\fa(N)\to\LH^n_\fa(M)\to\LH^n_\fa(M/N)\to 0$, if
$\LH^n_\fa(N)\neq 0$, then there exists
$\fp\in\Att(\LH^n_\fa(N))$. Thus, by Theorem 1.1, $\fp\in\Ass (N)$
and $\cd(\fa, R/\fp)=n$. As $\fp\in\Ass (M)$, by Theorem 1.1,
$\fp\in\Att(\LH^n_\fa(M))$, which contradicts with $\fp\in\Ass
(N)$. Therefore $\LH^n_\fa(N)=0$, and hence
$\LH^n_\fa(M)\cong\LH^n_\fa(M/N)$. Finally, as $\fa\subseteq\fm$,
by Theorem 1.5, there exists an epimorphism
$\LH^n_\fm(M/N)\to\LH^n_\fa(M/N)$, and the proof is complete.\hfill $\Box$\\

In [C. Corollary 1.7], Call shows that the set \{ $\LH
_{\fa}^{\dim (R)}(R): \fa \subseteq R \}$ has (up to isomorphism)
only finitely many modules. In the following result we give a
generalization of Call's result.

\subsection{Proposition}
{\it If
 ${\Att}({\LH}_{\fa}^{n}(M))= {\Att}({\LH}_{\fb}^{n}(M))$, then
  ${\LH}_{\fa}^{n}(M)\cong {\LH}_{\fb}^{n}(M)$. In addition the
   number of non--isomorphic top local cohomology modules
 ${\LH}_{\fa}^{n}(M)$ is not greater than $2^{|\Assh (M)|}-1$.}\\

{\it Proof}.  Since ${\LH}_{\fa}^{n}(M)$ and ${\LH}_{\fb}^{n}(M)$
have the structures of $ \widehat{R}$--modules, we may assume that
$R$ is complete. We take $N$ to be a submodule of $M$ such that
${\Ass}(M/N) ={\Att}({\LH}_{\fa}^{n}(M))$, ${\Ass}(N) =
{\Ass}(M)\setminus {\Att}({\LH}_{\fa}^{n}(M))$, and consider the
exact sequences
 $${\LH}_{\fa}^{n}(N)\longrightarrow
{\LH}_{\fa}^{n}(M)\longrightarrow
 {\LH}_{\fa}^{n}(M/N)\longrightarrow 0,$$
    $${\LH}_{\fb}^{n}(N)\longrightarrow
{\LH}_{\fb}^{n}(M)\longrightarrow
 {\LH}_{\fb}^{n}(M/N)\longrightarrow 0.$$

  \noindent So, as in the proof of Theorem 1.3, we have ${\LH}_{\fa}^{n}(M)\cong
 {\LH}_{\fm}^{n}(M/N)\cong {\LH}_{\fb}^{n}(M)$. \hfill $\Box$\\

 The following result is now clear for finite module $M$ with
 $|\Assh (M)|=1$.\\

 \subsection{Corollary}
  {\it Assume that ${\Assh}(M)$ is a singleton set. For any ideal $\fa$ if $\cd(\fa, M)=n$, then
  the non--zero local cohomology module ${\LH}_{\fa}^{n}(M)$ is isomorphic to ${\LH}_{\fm}^{n}(M)$.}\\

 {\it Proof}. The assertion follows from the facts that ${\LH}_{\fm}^{n}(M) \neq
 0$ and the number of non-isomorphic top local cohomology modules ${\LH}_{\fa}^{n}(M)$ is not greater than 1.\hfill
 $\Box$\\

 Now we are ready to present a reduced secondary representation of
 the top local cohomology modules of $M$ with respect to an ideal.
 Let $\Att(\LH^n_\fa(M))=\{\fp_1,\fp_2,\cdots ,\fp_r\}$. For each
 $1\le i\le r$, there exists a submodule $L_i$ of $M$ such that
 $\Ass (L_i)=\{\fp_i\}$ and $\Ass (M/L_i)=\Ass (M)\backslash\{\fp_i\}$.
 Consider the following exact sequences
 $$\LH^n_\fa(L_i)\stackrel{\varphi_i}{\to}\LH^n_\fa(M)\to\LH^n_\fa(M/L_i)\to
 0.$$
 Note that $\varphi_i(\LH^n_\fa(L_i))\neq 0$, otherwise
 $\LH^n_\fa(M)\cong\LH^n_\fa(M/L_i)$ but $\fp_i$ does not belong
 to $\Att(\LH^n_\fa(M/L_i))$.

 \subsection{Theorem} {\it With the above notations the following is a
 reduced secondary representation of $\LH^n_\fa(M)$,
 $$\LH^n_\fa(M)=\sum_{i=1}^r\varphi_i(\LH^n_\fa(L_i)).$$}\\

{\it Proof.} Let $1\le i\le r$ and $\bar{}:R\to R/\ann M$ be the
natural map. Since $\LH^n_\fa(L_i)\cong\LH^n_{\fa \bar{R}}(L_i)$,
$\Att_{\bar{R}}(\LH^n_{\fa\bar{R}}(L_i))=\{\bar{\fp_i}\}$ and
$\bar{\fp_i}$ is minimal in $\Spec \bar{R}$ we have that
$\LH^n_{\fa\bar{R}}(L_i)$ is a secondary $\bar{R}$-module, cf. [R,
Lemma 2]. Thus $\LH^n_\fa(L_i)$ is a secondary $R$-module too.
Therefore $\varphi_i(\LH^n_\fa(L_i))$ is also a $\fp_i$-secondary
$R$-module. Now by using the argument as in the proof of [DY,
Theorem
B] the assertion holds.\hfill$\Box$\\

\section{Last non--zero local cohomology modules}

This section is devoted to study the last non--zero local
cohomology modules $\LH _{\fa}^{\cd (\fa, M)}(M)$. We first prove
that $\cd (\fa, M)$ is equal to the infimum of integers $i\geq 0$
such that $\LH _{\fa}^{i+1}(R/\fp)$ is zero for all $\fp\in \Supp
(M)$ (see Theorem 2.1). We do not know much about the module $\LH
_{\fa}^{\cd (\fa, M)}(M)$, and the only information that we could
find is to specify a certain subset of the set of its attached
primes (see Theorem 2.2).

In [H, Proposition 2.3], Hartshorne showed that for an integer
$r$, if $\LH _{\fa}^{r}(N)= 0$ for all $R$--modules $N$, then the
corresponding condition holds for all $j\geq r$. In the following
Theorem we give a result related to the Hartshorne's one.

\subsection{Theorem}
{\it Let $r\geq 0$ be a given integer such that $\LH
_{\fa}^{r}(R/\fp)=0$ for all $\fp\in \Supp
(M)$. Then\\
\indent (a) $\LH _{\fa}^{r}(N)= 0$ for all finite $R$--modules $N$
with $\Supp (N) \subseteq \Supp (M)$. In particular $\LH
_{\fa}^{r}(M)= 0$.\\
\indent (b) $\cd (\fa, R/\fp)< r$ for all $\fp\in \Supp (M)$. In
particular $\cd (\fa, M)< r$.}\\

{\it Proof}. (a). It is enough to show that $\LH _{\fa}^{r}(M)=
0$. There is a prime filtration $0= M_{0}\subseteq M_{1}\subseteq
\cdots \subseteq M_{s}= M$ of submodules of $M$ such that, for
each $1\leq j\leq s$, $M_{j}/M_{j-1}\cong R/\fp _{j}$, where $\fp
_{j}\in \Supp (M)$. From the exact sequences $$\LH
_{\fa}^{r}(M_{j-1})\longrightarrow \LH
_{\fa}^{r}(M_{j})\longrightarrow \LH _{\fa}^{r}(R/\fp _{j})$$ we
eventually get $\LH _{\fa}^{r}(M)= 0$.\\
\indent (b). By (a), $\LH _{\fa}^{r}(R/\fb)= 0$ for all ideals
$\fb$ of $R$ with $V(\fb)\subseteq \Supp (M)$. We first show that
$\LH _{\fa}^{r+1}(R/\fp)= 0$ for all $\fp\in \Supp (M)$. Assume
that $\LH _{\fa}^{r+1}(R/\fp)\neq 0$ for some $\fp\in \Supp (M)$,
so that $\fa\nsubseteq \fp$. Choose $0\neq w\in \LH
_{\fa}^{r+1}(R/\fp)$. Since $w$ is annihilated by some power of
$\fa$, there exists $x\in \fa\setminus \fp$ such that $xw= 0$. Now
consider the exact sequence $0\longrightarrow
R/\fp\overset{x}\longrightarrow R/\fp\longrightarrow
R/(\fp+xR)\longrightarrow 0$ which induces the exact sequence
$$\LH _{\fa}^{r}(R/(\fp+xR))\longrightarrow \LH
_{\fa}^{r+1}(R/\fp)\overset{x}\longrightarrow \LH
_{\fa}^{r+1}(R/\fp).$$ As $V(\fp+ xR)\subseteq \Supp (M)$, by (a),
$\LH _{\fa}^{r}(R/(\fp+xR))= 0$. Thus the map $\LH
_{\fa}^{r+1}(R/\fp)\overset{x}\longrightarrow \LH
_{\fa}^{r+1}(R/\fp)$ is injective. As $xw= 0$, we get $w= 0$ which
is a contradiction. Therefore $\LH _{\fa}^{r+1}(R/\fp)= 0$ for all
$\fp\in \Supp (M)$. By (a), $\LH _{\fa}^{r+1}(N)= 0$ for all
finite $R$--modules $N$ with $\Supp (N) \subseteq \Supp (M)$. We
proceed by induction and the proof is complete. \hfill $\Box$

\subsection{Corollary}
{\it For a finite $R$-module $M$, we have:}\\
  \begin{center} $\cd(\fa, M)= \inf \{r\geq 0 :
  \LH _{\fa}^{r}(R/\fp)= 0$ for all
$\fp\in \Supp (M) \}-1$ . \end{center} \hfill $\Box$\\

\subsection{Corollary}
{\it For each integer $i$ with $0\leq i\leq \cd (\fa, M)$, there
exists $\fp\in \Supp (M)$ with $\LH _{\fa}^{i}(R/\fp)\neq 0$}.
\hfill $\Box$\\

Our final result specifies a certain subset of $\Att (\LH
_{\fa}^{\cd (\fa, M)}(M))$.\\

 \subsection{Theorem}
 {\it The set} $T:= \{ \fp\in \Ass (M) : \cd (\fa, R/\fp)= \dim (R/\fp)= \cd (\fa,
 M)\}$ {\it is a subset of} $\Att (\LH_{\fa}^{\cd (\fa,
 M)}(M))$. {\it In particular, we have equality if $\cd(\fa, M)=\dim
 (M)$.}\\

 {\it Proof}. Assume that the set $S:=\{ \fp\in \Ass (M) : \dim (R/\fp)= \cd (\fa,
 M)\}$ is not empty. There exists a submodule $N$ of $M$ with
 $\Ass (N)= \Ass (M)\setminus S$ and $\Ass(M/N)= S$. As $\Supp
 (N)\subseteq\Supp (M)$, we have the exact sequence\\
 $$\LH_{\fa}^{\cd (\fa, M)}(N)\longrightarrow \LH_{\fa}^{\cd (\fa,
 M)}(M)\longrightarrow \LH_{\fa}^{\cd (\fa, M)}(M/N)\longrightarrow 0$$
 (see [DNT, Theorem 2.2]). This shows that $\Att (\LH_{\fa}^{\cd (\fa,
 M)}(M/N))\subseteq \Att (\LH_{\fa}^{\cd (\fa,
 M)}(M))$. Therefore it is sufficient to show that $T= \Att (\LH_{\fa}^{\cd (\fa,
 M)}(M/N))$. We have $\dim (R/\fp)= \cd (\fa,
 M)$ for all $\fp\in \Ass (M/N)$. Therefore $\cd (\fa, M)= \dim
 (M/N)$,and hence $\LH_{\fa}^{\cd (\fa,
 M)}(M/N)= \LH_{\fa}^{\dim (M/N)}(M/N)$ is an Artinian $R$--module
 and, by Theorem 1.1, $ \Att (\LH_{\fa}^{\cd (\fa,
 M)}(M/N))= T$.

 For the special case see Theorem 1.1.\hfill $\Box$\\

{\bf \small {The Authors' Addresses}} \\[0.3 cm]
{\footnotesize Mohammad T. Dibaei, Mosaheb Institute of
Mathematics, Teacher Traning University, Tehran, Iran.\\
E-mail address: {\tt dibaeimt@ipm.ir} \\[0.2cm]
Siamak Yassemi, School of Mathematics, Institute for Studies in
Theoretical Physics and Mathematics, Tehran, Iran, and\\
Department of Mathematics, University of Tehran, Tehran, Iran.\\
E-mail address: {\tt yassemi@ipm.ir}}

\end{document}